\documentclass[11pt]{article}
\usepackage{amsmath}
\usepackage{mathrsfs}
\usepackage{amsfonts}

\usepackage{amsfonts, amsmath, amssymb, amssymb}
\usepackage{amssymb,amsfonts,amsmath,
latexsym, epsfig,cite, psfrag,eepic,colordvi}
\usepackage{amscd,graphics}
\usepackage{lineno}
\newcommand{\qed}{\hfill $\Box $}
\newcommand{\pf}{\noindent {\bf Proof.} }

\newtheorem{theorem}{Theorem}[section]
\newtheorem{lemma}[theorem]{Lemma}

\begin{document}

\title{An Alternative Proof of the $H$-Factor Theorem
\thanks{This work is supported in part by Natural Sciences and Engineering Research Council of Canada.}}

\author{Hongliang Lu\,\textsuperscript{a}, \
Qinglin Yu\,\textsuperscript{b}\thanks{Corresponding email: yu@tru.ca (Q. Yu)}
\\ {\small \textsuperscript{a}Department of Mathematics}
\\ {\small Xi'an Jiaotong University, Xi'an,  China}
\\ {\small \textsuperscript{b}Department of Mathematics and Statistics}
\\ {\small Thompson Rivers University, Kamloops, BC, Canada}
%\\ {\small \textsuperscript{c}School of Mathematics, Shandong University, Jinan,
%China}
}

\date{}

\maketitle

\date{}

\maketitle

\begin{abstract}
Let $H: V(G) \rightarrow 2^{\mathbb{N}}$ be a set mapping for a
graph $G$. Given a spanning subgraph $F$ of $G$, $F$ is called a
{\it general factor} or an $H$-{\it factor} of $G$ if $d_{F}(x)\in
H(x)$ for every vertex $x\in V(G)$.  $H$-factor problems are, in
general, $NP$-complete problems and imply many well-known factor
problems (e.g., perfect matchings, $f$-factor problems and $(g,
f)$-factor problems) as special cases. Lov\'asz [The factorization
of graphs (II),
  Acta Math. Hungar.,
23 (1972), 223--246] gave a structure description and obtained a
deficiency formula for $H$-optimal subgraphs. In this note, we use a
generalized alternating path method to give a structural
characterization and provide an alternative and shorter proof of
Lov\'asz's deficiency formula.
\end{abstract}

\section{Introduction}

In this paper, we consider finite undirected graphs without loops and multiple edges. For a graph
$G = (V, E)$, the degree of $x$ in $G$ is denoted by $d_G(x)$, and the set of vertices adjacent to
$x$ in $G$ is denoted by $N_G(x)$. For $S\subseteq V(G)$, the subgraph of $G$ induced by $S$ is
denoted by $G[S]$ and $G-S= G[V(G)-S]$. For vertex subsets $S$ and $T$, $E_G(S, T)$ is the set of
edges between $S$ and $T$ in $G$.  We use $\omega(G)$ for the number of connected components in
$G$. Notations and terminologies not defined here may be found in \cite{lov3}.

For a given graph $G$, we associate an integer set $H(x)$ with each
vertex $x\in V(G)$ (i.e., $H$ is a set mapping from $V(G)$ to
$2^{\mathbb{N}}$). Given a spanning subgraph $F$ of $G$, $F$ is a
{\it general factor} or an $H$-{\it factor} of $G$ if $d_{F}(x)\in
H(x)$ for every vertex $x\in V(G)$. By specifying $H(x)$ to be an
interval or a special set, an $H$-factor becomes an $f$-factor, an
$[a, b]$-factor or a $(g, f)$-factor, respectively. For a general
mapping $H$, the decision problem of determining whether a graph has
an $H$-factor is known to be $NP$-complete. In fact, when $H(x)$
contains a ``gap" with more than one element, $H$-factor problem is
an $NP$-complete problem. Interestingly, Lov\'asz \cite{Lov73}
showed that the 3-edge-colorability problem is reducible to
$H$-factor problem with $H(x) = \{1\}$ or $\{0, 3\}$ for each $x \in
V$. Furthermore, he proved that Four-Colors Problem is equivalent to
the existence of such a special $H$-factor. So it is reasonable to
conclude that finding a characterization for $H$-factors in general
is a challenging problem and so it is natural to turn our attention
to $H$-factor problems in which $H(x)$ contains only one-element
gaps. Furthermore, Lov\'asz also conjectured that the general factor
problem with one-element gaps  could be solved in polynomial time
and Cornu\'ejols \cite{Cor88} proved the conjecture.

Assume that $H$ satisfies the property:
\begin{align*}
(*)\ \ \ \  \mbox{if}\ i\not\in H(x),\ \mbox{then}\ i+1\in H(x),\
\mbox{ for}\ mH(x)\leq i\leq MH(x),
\end{align*}
where $mH(x) = \min \{ r \ | \ r \in H(x)\}$ and $MH(x) = \max \{ r
\ | \ r \in H(x)\}$. Let $MH(S)=\sum_{u\in S}MH(u)$,
$mH(S)=\sum_{v\in S}mH(v)$ and $H\pm c = \{i \pm c \ | \ i \in H
\}$. Lov\'asz \cite{lov1} obtained a sufficient and necessary
condition for the existence of $H$-factors with the properties $(*)$
and a deficiency formula for $H$-optimal subgraphs.  In this paper,
we use the traditional technique -- alternating path -- which has
dealt effectively with other factor problems to prove Lov\'asz's
deficiency formula. However, we need to modify the usual alternating
paths to {\it
changeable trails} to handle the more complicated structures in this case. \\

Let $H$ be a set function satisfying the property (*) and $F$ any
{\it spanning} subgraph of $G$. For a vertex $x \in V(G)$, if
\begin{equation}\label{feasible}
d_F(x) \in H(x),
\end{equation}
then vertex $x$ is called {\it feasible}. So a subgraph $F$ is an $H$-factor if and only if every
vertex is feasible. Given a spanning subgraph $F$ and a subset $S \subseteq V(G)$, the {\it
 deficiency of subgraph} $G[S]$ in $F$ is defined as
$$def_H(F; S)=\sum_{x\in S}dist(d_{F}(x),H(x)),$$
where $dist(d_{F}(x),H(x))$ is the distance of $d_F(x)$ from the set
$H(x)$. For convenience, we use $def(F; S)$ in short. So $def(F; x)
= dist(d_F(x), H(x))$ is the deficiency of vertex $x$ in $F$.
We can
measure $F$'s ``deviation" from condition (\ref{feasible}) by
defining the \emph{deficiency} of $F$ with respect to $H$ as
\begin{align*}
def_H[F]=\sum_{x\in V(F)}dist(d_{F}(x),H(x)).
\end{align*}
The {\it total deficiency} of $G$ with respect to $H$ is
$$def_H(G) = \min \{def_H[F] \ | \  \mbox{$F$ is a spanning subgraph of $G$} \}.$$
Clearly, $def_H(G) = 0$ if and only if there exists an $H$-factor. A
subgraph $F$ is called $H$-{\it optimal}, if $def_H[F] = def_H(G)$.
Of course, any $H$-factor is $H$-optimal.

Let $I_H(x) = \{ d_F(x) \ | \ \mbox{$F$ is any $H$-optimal subgraph}
\}$. Lov\'asz \cite{lov1} studied the structure of $H$-factors in
graphs by introducing a Gallai-Edmonds type of partition for $V(G)$
as follows:

$$\begin{array}{cll}
C_H(G) & = & \{x \ | \ I_H(x) \subseteq H(x)\}, \\
A_H(G) & = &\{x\ |\ \min I_H(x) \geq MH(x)\},\\
B_H(G) & = & \{x\ |\ \max I_H(x) \leq mH(x)\},\\
D_H(G) & = & V(G) - A_H(G) - B_H(G) - C_H(G).
\end{array}$$

Based on this canonical partition, Lov\'asz obtained a sufficient
and necessary conditions of $H$-factors with property (*) and as
well as the deficiency formula for $H$-optimal subgraphs. In this
paper, we give an alternative description of the partition $(A, B,
C, D)$ by deploying changeable trails and therefore provide a new
proof of the deficiency formula for $H$-optimal subgraphs. Our
approach is as follows:

Suppose that $G$ does not have $H$-factors. Choose a spanning
subgraph $F$ of $G$ such that for all $v\in V(G) , d_{F}(v)\leq
MH(v)$ and the deficiency is minimized over all such choices.
Moreover, we choose $F$ such that the $E(F)$ is {\it minimal}.
Necessarily, there is a vertex $v\in V$ such that $d_{F}(v)\not\in
H(v)$, so the deficiency of $F$ is positive.
  Set
\begin{align*}
B_{0}=\{x\in V(G)\ |\ d_{F}(x)\not\in H(x)\}.
\end{align*}
Since $E(F)$ is minimal and $H$ satisfies (*), we have
\begin{align*}
B_{0}=\{x\in V(G)\ |\ d_{F}(x)<mH(x)\}.
\end{align*}

A trail $P=v_{0}v_{1}\dots v_{k}$ is called a \emph{changeable
trail} if it satisfies the following condition:
\begin{itemize}
\item[(a)]$v_{0}\in B_{0}$, and $v_{0}v_{1}\notin E(F)$;

\item[(b)]$def(F; x)=def(F\bigtriangleup P; x)=0$, for every $x\in
V(P)-v_{0}-v_{k}$;

\item[(c)] if $v_{0}= v_i \neq v_{k}$, then $def(F; v_0) > def(F\bigtriangleup P; v_0)$;

\item[(d)] for all $l\leq k$, sub-trail $P'=v_{0}v_{1}\dots v_{l}$ satisfies conditions (a)-(c) as well.
\end{itemize}
A changeable trail $P$ is {\it odd} if the last edge doesn't belong
to $F$; otherwise, $P$ is {\it even}. Moreover, the trails of length
zero are considered as even changeable trails.

For a given graph $G$, we define $D(G)$ to be a vertex set consisting
of three types of vertices as follows:
\begin{itemize}
\item[(i)]$\{v\ |\ \exists\ \mbox{both of an even changeable  trail and an odd
changeable trail} \ \mbox{ from}\ B_{0} \ \mbox{to}  \
 v\}$;

\item[(ii)]$\{v\ |\ mH(v)<d_{F}(v)\leq MH(v)\ \mbox{and}\ \exists\
 \mbox{an even changeable  trail}\ \mbox{from}\ B_{0}\ \mbox{to}  \
 v\}$;

\item[(iii)]$\{v\ |\ mH(v)\leq d_{F}(v)<MH(v)\ \mbox{and}\
 \exists\mbox{ an odd changeable  trail}\ \mbox{from}\ B_{0}\ \mbox{to}  \
 v\}$.
\end{itemize}
%\textbf{(In fact, a vertex $v\in D$ can be type (i), (ii), (iii)).}
The sets $A(G)$ and $B(G)$ are defined as follows:
\begin{align*}
B(G)=\{v\ |\ \exists\mbox{ an even changeable trail ending at} \
v\}-D,\\
A(G)=\{v\ |\ \exists\mbox{ an odd changeable trail ending at}\
v\}-D,
\end{align*}
and $C(G)=V(G)-A(G)-B(G)-D(G)$. We abbreviate $D(G), A(G), B(G)$ and
$C(G)$ by $D, A, B$ and $C$, respectively.

If $v\in B$, then $d_{F}(v)\leq mH(v)$.  Otherwise, as $v \not\in
D$, we can swap edges in $F$ along an even changeable trail ending
at $v$ and thus decrease the deficiency. Similarly, if $v\in A$,
then $d_{F}(v)=MH(v)$.  Otherwise, as $v\not\in D$, we can likewise
decrease the deficiency by swapping edges in $F$ along an odd
changeable trail ending at $v$. By the definitions, clearly $A$,
$B$, $C$ and $D$ are a partition of $V(G)$. We call a changeable
trail $P$ an \emph{augmenting changeable trail} if $def(F\triangle
P; G)<def(F; G)$. Following the above discussion, when $H(v)$ is an
integer interval with more than an element, then $v\notin D$.

\section{Main Theorem}

In the following lemmas, we always assume that $G$ has no
$H$-factors and $F$ is an $H$-optimal subgraph with minimal $E(F)$.
Let $\tau = \omega(G[D])$ and $D_{1}, \dots, D_{\tau}$ be the
components of the subgraph of $G$ induced by $D$.

\begin{lemma}\label{lem1}
An $H$-optimal subgraph $F$ does not contain an augmenting changeable trail.
\end{lemma}

\begin{lemma}\label{prop1}
$def(F; D_{j})\leq 1$ for $j=1,\ldots,\tau$.
\end{lemma}
\vspace{-3mm}\pf Suppose, to the contrary, that $def(F; D_{i}) > 1$. Let $v_0\in {D_{i}}$ and
$def(F; v_0)\geq 1$. Since $E(F)$ is minimal, so $d_{F}(v_0)<mH(v_0)$. Hence $v_0$ is of type (i)
and there exists an odd changeable trail $P$ from a vertex $x$ of $B_{0}$ to $v_0$. Then $x=v_0$.
Otherwise, $def(F; G)>def(F\Delta P; G)$, a contradiction since $F$ is $H$-optimal. Furthermore, if
$def(F; v_0)\geq 2$, then $def(F; G)>def(F\Delta P; G)$, a contradiction again. So we have $def(F;
v_0)=1$ and $def(F; u)\leq 1$ for any $u\in D_{i}-v_0$. Moreover, $d_{F}(v_0)+1\in H(v_{0})$ and
$d_{F}(v_0)+2\notin H(v_{0})$.

We define $D_{i}^{1}$ to be a vertex set consisting of three types of
vertices as follows:
\begin{itemize}
\item[(1)]$\{w\in D_{i}\ |\ \exists\mbox{ an even changeable  trail and an
odd changeable trail} \ \mbox{from}\ v_0 \ \mbox{to}  \
 w\}$;

\item[(2)]$\{w\in D_{i}\ |\ mH(w)<d_{F}(w)\leq MH(w)\
 \mbox{and}\ \exists\mbox{ an even changeable  trail}\ \mbox{from}\ v_0\ \mbox{to}  \
 w \} $;

\item[(3)]$\{w\in D_{i}\ |\ mH(w)\leq d_{F}(w)<MH(w)\
 \mbox{and}\ \exists\mbox{ an odd changeable  trail}\ \mbox{from}\ v_0\ \mbox{to}  \
 w\}$.
\end{itemize}
%%%%%%%%%%%%%%%%%有的顶点可能同时满足三种types
Now we choose a maximal subset $D_{i}^{2}$  of $D_{i}^{1}$ such that
$P \subseteq D_{i}^{2}$ and  the trails, which are  of type (1), (2)
or (3), belongs to $D_{i}^{2}$.

\medskip
{\it Claim.} $D_{i}^{2}=D_{i}$.
\medskip

Otherwise, since $D_{i}$ is connected, there exists an edge $xy \in E(G)$ such that $x\in
D_{i}-V(D_{i}^{2})$ and $y\in V(D_{i}^{2})$. We consider $xy\in E(F)$ (or $xy \not\in E(F)$).

Then there exists an even (resp. odd) changeable trail $P_{1}$ from
$v_0$ to $x$, where $xy\in P_{1}$ and $V(P_{1})-x\subseteq
V(D_{i}^{2})$.
%There exists an even (resp. odd) changeable trail from $v_0$ to $x$. %%%%这句就是感觉重复了，所以改掉
 So $x$ is type (i) or type
(ii) (resp. type (i) or type (iii)). Since $x\notin D_{i}^2 $, $x$
can only be of type (i).  %Since $x\in D_{i}$ and $x$ is of type (i), %%%%这个跟上面一样
So there exists an odd (resp. even) changeable trail $P_2$ from a
vertex $t$ of $B_{0}$ to $x$. Thus $t\neq v_0$; otherwise, we have
$V(P_{1}\cup P_{2})\subseteq D_{i}^{2}$, a contradiction to the
maximality of $D_{i}^{2}$. If $E(P_{1})\cap E(P_{2})=\emptyset$,
then $def(F; G)>def(F\bigtriangleup (P_{1}\cup P_{2}); G)$, a
contradiction since $F$ is $H$-optimal. Let $z\in P_{2}$ be the
first vertex which also belongs to $D_{i}^{2}$ and denote the
subtrail of $P_2$ from $t$ to $z$ by $P_{3}$. If $z$ is of type (1),
by the definition, there exist both an odd changeable trail $P_{4}$
from $v$ to $z$ and an even changeable trail $P_{5}$ from $v$ to $z$
such that $V(P_{4}\cup P_{5})\subseteq V(D_{i}^{2})$. Thus either
$P_{4}\cup P_{3}$ or $P_{5}\cup P_{3}$ is an augmenting trail, a
contradiction to Lemma \ref{lem1}. If $z$ is type (2) or type (3),
the argument is similar. We complete the claim.

\medskip

Let $u\in V(D_{i})-v_0$ and $def(F; u)=1$. Since $u$ is not type (2),  there exists an odd
changeable trail $P_{6}$ from $v_0$ to $u$. We have $def(F; G)>def(F\bigtriangleup P_{6}; G)$, a
contradiction since $F$ is $H$-optimal. \qed

Using the above lemma, we have the following result.

\begin{lemma}\label{prop2}
 For $i=1,\ldots,\tau$, if $def(F; D_{i})=1$, then
\begin{itemize}
\item[$(a)$] $E_G(D_{i},B)\subseteq E(F)$;

\item[$(b)$] $E_G(D_{i},A)\cap E(F)=\emptyset$.
\end{itemize}
\end{lemma}
\vspace{-3mm} \pf Let $def(F; r)=1$, where $r\in V(D_{i})$. Suppose
the lemma does not hold.

To show (a), let $uv\not\in E(F)$, where $u\in V(D_{i})$ and $ v\in V(B)$. If $u$ is of type (i) or
type (ii), from the proof of Lemma \ref{prop1}, then there exists an even changeable trail
$P\subseteq G[D_{i}]$ from $r$ to $u$. Hence $P\cup uv$ be an odd changeable trail from $r$ to $v$,
a contradiction to $ v\in B$. If $u$ is of type (iii),  then there exists an odd changeable trail
$P\subseteq G[D_{i}]$ from $r$ to $u$. Since $F$ is $H$-optimal and $H$ has the property (*), so
$d_{F}(u)\in H(u)$, $d_{F}(u)+1\notin H(u)$ and $d_{F}(u)+2\in H(u)$. Hence $P\cup uv$ is an odd
changeable trail from $r$ to $v$, a contradiction to $ v\in B$ again.

%%%%%%%%%%%%%%后面这段是加上的，补充完整。由于不好叙述，上面的之前
%%%%%%%%%%%%%%%写的两个cases给合到一块了

Next we consider (b). Let $uv\in F$, where $u\in D_{i}$ and $v\in A$. If $u$ is of type (i) or type
(iii), from the proof of Lemma \ref{prop1}, then there exists an odd changeable trail $P\subseteq
G[D_{i}]$ from $r$ to $u$. Then $P\cup uv$ be an even changeable trail from $r$ to $v$, a
contradiction to $ v\in A$. If $u$ is of type (ii), then there exists an even changeable trail
$P\subseteq G[D_{i}]$ from $r$ to $u$. Since $F$ is $H$-optimal and $H$ has the property (*), so
$d_{F}(u)\in H(u)$, $d_{F}(u)-1\notin H(u)$ and $d_{F}(u)-2\in H(u)$. Hence $P\cup uv$ is an even
changeable trail from $r$ to $v$, a contradiction to $ v\in A$ again.
 \qed

From the definition of partition $(A, B, C, D)$, it is not hard to
see the next lemma.
\begin{lemma}
$E_G(B,C\cup B)\subseteq E(F)$, $E_{G}(A,A\cup C)\cap
E(F)=\emptyset$ and $E_G(D,C)=\emptyset$.
\end{lemma}

%\begin{lemma}\label{prop4}
%$F$ misses at most one edge of $E_G(D_{i}, B)$ and contains at most
%one edge of $E_G(D_{i}, A)$. Moreover, if $F$ misses exactly one
%edge of $E_G(D_{i}, B)$, then $E_G(D_{i},A)\cap E(F)=\emptyset$; if
%$F$ contains  exactly one edge of $E_G(D_{i}, A)$, then
%$E_G(D_{i},B)\subseteq E(F)$.
%\end{lemma}

\begin{lemma}\label{prop4}
\begin{itemize}
\item[$(a)$] $F$ misses at most one edge of $E_G(D_{i}, B)$. Moreover, if $F$ misses  one
edge of $E_G(D_{i}, B)$, then $E_G(D_{i},A)\cap E(F)=\emptyset$;

\item[$(b)$] $F$ contains at most
one edge of $E_G(D_{i}, A)$. Moreover, if $F$ contains   one edge of $E_G(D_{i}, A)$, then
$E_G(D_{i},B)\subseteq E(F)$.
\end{itemize}
\end{lemma}
\vspace{-3mm} \pf By Lemma \ref{prop2}, we may assume $def(F;
D_{i})=0$. Let $u\in V(D_{i})$, by the definition of $D$, there
exists a changeable trail $P$ from a vertex $x$ of $B_{0}$ to $u$.
Denote the first vertex in $P$ belonging to $D_{i}$ by $y$, and the
sub-trail of $P$ from $x$ to $y$ by $P_{1}$. Let $y_{1}y\in
E(P_{1})$, where $y_1\notin D_i$. Without loss of generality, assume
that $P_{1}$ is an odd changeable trail (when $P_{1}$ is an even
changeable trail, the proof is similar). Since $P_{1}$ is a
changeable trail, so $y_{1}\in B$ and $y_{1}y\notin F$. Because
$y\in D$, $y$ is of type (i) or type (iii). We define the subset
$D_{i}^{1}\subseteq D_{i}$ which consists of the following vertices:
\begin{itemize}
\item[$(1)$]$\{w\in D_{i}\ |\ \exists\mbox{ an even changeable  trail and an
odd changeable trail along $P_{1}$} \ \mbox{ from}\ x \ \mbox{to}  \
 w\}$;

\item[$(2)$]$\{w\in D_{i}\ |\ mH(w)<d_{F}(w)\leq MH(w)\
 \mbox{and}\ \exists\mbox{ an even changeable  trail along $P_{1}$}\ \mbox{from}\ x\ \mbox{to}  \
 w \} $;

\item[$(3)$]$\{w\in D_{i}\ |\ mH(w)\leq d_{F}(w)<MH(w)\
 \mbox{and}\   \exists\mbox{  an odd changeable  trail along $P_{1}$}\ \mbox{from}\ x\ \mbox{to}  \
 w\}$.
\end{itemize}

Now we choose a maximal subset $D_{i}^{2}$  of $D_{i}^{1}$ such that
the trails, which are of type (1), (2) or (3), except $V(P_{1})-y$,
belongs to $D_{i}^{2}$.

\medskip
{\it Claim 1.} $D_{i}^{1}=D_{i}=D_{i}^{2}$.

Suppose that $D_{i}\neq D_{i}^{2}$. Let $v_{1}v_{2}\in E(G)$, where
$v_{1}\in D_{i}^{2}$ and $v_{2}\in D_{i}-D_{i}^{2}$. Firstly, we
show that $D_{i}^{2}\neq \emptyset$. If $y$ is type (iii), then
$y\in D_{i}^{2}$. If $y$ is type (i) or type (ii), then there exists
an even changeable trail $R_{1}$ from a vertex $w$ of $B_{0}$ to
$y$. We have $yy_{1}\in E(R_{1})$; otherwise $R_{1}\cup yy_{1}$ is
an odd changeable trail from $w$ to $y_{1}$, contradicting to
$y_{1}\in B$. Hence, we may assume $w=x$ and $P_{1}$ is a subtrail
of $R_{1}$. So $V(R_{1})-(V(P_{1})-y)\subseteq D_{i}^{2}$ and
$D_{i}^{2}\neq \emptyset$.

%%%%%%%
%%Lu 接下来我标出了几处黑体，是想说明两种情况基本上是相似的，但
%%%并没完全的标出来。接下来这段跟lemma 2.2中那个claim证明就非常的相似点了
%%%%%%

We consider $v_{1}v_{2}\in E(F)$. Then there exists an even
changeable trail $R_{2}$ from $x$ to $v_{2}$ such that
$V(R_{2})-(V(P_{1})-y)-v_{2}\subseteq D_{i}^{2}$. If $v_{2}$ is type
(ii), by the definition of $D_{i}^2$, then we have $v_{2}\in
D_{i}^2$, contradicting to the maximality of $D_{i}^{2}$. If $v_{2}$
is type (i) or type (iii), then there exists an odd changeable trail
$R_{3}$ from a vertex $w_{2}$ of $B_{0}$ to $v_{2}$. Next we show
that $yy_{1}\in R_{3}$. If $V(R_{3})\cap D_{i}^2\neq \emptyset$, let
$z$ be first vertex in $R_{3}$ belonging $D_{i}^2$; else let
$z=v_{2}$. Without loss of generality, we suppose that the subtrail
$R_{4}$ from $w_{2}$ to $z$ along $R_{3}$ is an odd changeable trail
and $z\in D_{i}^2$. If $z\in D_{i}^2$ is type (1) or (2), then there
is an even changeable trail, say $R_{5}$, from $x$ to $z$ along
$P_{1}$ such that $V(R_{5})-(V(P_{1})-y)\subseteq D_{i}^2$. Let
$R_{6}$ is a subtrail from $y_{1}$ to $z$ along $R_{5}$. Then
$R_{4}\cup R_{6}$ is an odd changeable trail from $w_{2}$ to
$y_{1}$, contradicting to $y_{1}\in B$. If $z\in D_{i}^2$ is type
(3), then $d_{F}(z)\in H(z)$, $d_{F}(z)+1\notin H(z)$, and
$d_{F}(z)+2\in H(z)$. Moreover, there is an odd changeable trail
$R_{7}$ along $P_{1}$ from $x$ to $z$. Let $R_{8}$ is a subtrail
from $y_{1}$ to $z$ along $R_{7}$. Then $R_{4}\cup R_{8}$ is an odd
changeable trail from $w_{2}$ to $y_{1}$, contradicting to $y_{1}\in
B$ again.  So $yy_{1} \in R_{3}$. Let $R_{9}$ be the subtrail from
$y$ to $v_{2}$ along $R_{3}$. Then we have $V(R_{9})\subseteq
D_{i}^2$, contradicting to the maximality of $D_{i}^2$. By the
symmetry of definition of $D_{i}$ and $D_{i}^2$,   for
$v_{1}v_{2}\notin E(F)$, the proof is similar. We complete the proof
of the claim.

Let $x_{3}y_{3}\in E(G)-yy_{1}$. We have the following two claims.

\medskip
{\it Claim 2.} If $x_{3}\in D_{i}$ and $y_{3}\in B$, then
$x_{3}y_{3}\in E(F)$.

Otherwise,  $x_{3}y_{3}\notin E(F)$. If $x_{3}$ is of type (1) or
(2), by the definition of set $D_{i}^2$ and $D_{i}=D_{i}^{2}$, then
there exists an even changeable trail $P_{10}$ from $x$ to $x_{3}$
such that $V(P_{10})-(V(P_{1})-y)\subseteq D_{i}$. Then $ P_{10}\cup
x_{3}y_{3}$ is an odd changeable trail from $x$ to $y_{3}$,
contradicting to $y_{3}\in B$.
%\medskip
%{\it Case 2.} $x_{3}$ is of type (2).
%
%By the definition, there exists an even changeable trail $P_{5}$
%from $x$ to $x_{3}$ such that $V(P_{5})-(V(P_{1})-y)\subseteq
%D_{i}$. Note that  $d_{F}(x_{3})\in H(x_{3})$. Since $F$ is
%$H$-optimal and $H$ has the property (*), so $d_{F}(x_{3})-1\notin
%H(x_{3})$ and $d_{F}(x_{3})-2\in H(x_{3})$. Then $ P_{5}\cup
%x_{3}y_{3}$ is an odd changeable trail from $x$ to $y_{3}$, a
%contradiction.
If $x_{3}$ is of type (3), then there exists an odd changeable trail
$P_{11}$ from $x$ to $x_{3}$ such that
$V(P_{11})-(V(P_{1})-y)\subseteq D_{i}$. Note that  $d_{F}(x_{3})\in
H(x_{3})$. Since $F$ is $H$-optimal, so $d_{F}(x_{3})+1\notin
H(x_{3})$ and $d_{F}(x_{3})+2 \in H(x_{3})$. Then $P_{11}\cup
x_{3}y_{3}$ is an odd changeable trail from $x$ to $y_{3}$,
contradicting to $y_{3}\in B$. We complete Claim 2.

\medskip
{\it Claim 3.} If $x_{3}\in D_{i}$ and $y_{3}\in A$, then
$x_{3}y_{3}\notin E(F)$.

Otherwise,  $x_{3}y_{3}\in E(F)$. If $x_{3}$ is of type (1) or (3),
then there exists an odd changeable trail $P_{12}$ from $x$ to
$x_{3}$ such that $V(P_{12})-(V(P_{1})-y)\subseteq D_{i}$. Then
$P_{12}\cup x_{3}y_{3}$ is an even changeable trail from $x$ to
$y_{3}$, contradicting to $y_{3}\in A$. If $x_{3}$ is of type (2),
then there exists an even changeable trail $P_{13}$ from $x$ to
$x_{3}$ such that $V(P_{13})-(V(P_{1})-y)\subseteq D_{i}$. Note that
$d_{F}(x_{3})\in H(x_{3})$. Since $F$ is $H$-optimal and $H$ has the
property (*), so $d_{F}(x_{3})-1\notin H(x_{3})$ and
$d_{F}(x_{3})-2\in H(x_{3})$. Then $ P_{13}\cup x_{3}y_{3}$ is an
even changeable trail from $x$ to $y_{3}$, contradicting to
$y_{3}\in A$ again.

We complete the proof. \qed

Now we present and prove deficiency formula for $H$-optimal
subgraphs. Recall that $\tau$ is the number of components in $G[D]$.

\begin{theorem}\label{thm1}
$def_{H}(G)=\tau+\sum_{v\in B}(mH(v)-d_{G-A}(v))-\sum_{v\in
A}MH(v)$.
\end{theorem}
\vspace{-2mm}\pf Let $\tau_{1}$ denote the number of components
$D_{i}$ of $G[D]$ which satisfies $def(F; D_{i})=1$. Let $\tau_{B}$
(or $\tau_{A}$) be the number of components $T$ of $G[D]$ such that
$F$ misses (or contains) one edge from $T$ to $B$ (or $A$). By
Lemmas \ref{prop2} and \ref{prop4}, we have
$\tau=\tau_{1}+\tau_{A}+\tau_{B}$. Note that $d_{F}(v)\leq mH(v)$
for all $v\in B$ and $d_{F}(v)=MH(v)$ for all $v\in A$. So
\begin{align*}
def_{H}(G)&=\tau_{1}+mH(B)-\sum_{v\in B}d_{F}(v)\\
 &=\tau_{1}+mH(B)-(\sum_{v\in
B}d_{G-A}(v)-\tau_{B})-(MH(A)-\tau_{A})\\
&=\tau+mH(B)-\sum_{v\in B}d_{G-A}(v)-MH(A).
\end{align*}  \qed

%A graph $G$ is called $H$-{\it critical} if $V(G)=D$.
Let $X, Y $ be two disjoint subsets  of $V(G)$.  Define the modified prescription of $H$ to be
\begin{align*}
H_{(X,Y)}(u)=H(u)-|E_G(u,Y)|\ \ \mbox{for $u\in V(G)-X-Y$}.
\end{align*}
Let $K$ be a component of $G-X-Y$. We defined $H_{(X,Y)|K}$ as follows:
\begin{align*}
H_{(X,Y)|K}(u)=H_{(X,Y)}(u)\ \ \mbox{for $u\in V(K)$}.
\end{align*}
%
%By above discussion of Lemmas \ref{prop1} and \ref{prop4},  we
%obtain the following result.
%\begin{theorem}
%Every component $D_{i}$ of $D$ is $H_{(A,B)|D_{i}}$-critical for
%$i=1,\ldots,\tau$.
%\end{theorem}

%For any pair of disjoint subsets $S,T\subseteq V(G)$ and  let
%$\tau_{H}(S,T)$ denote the  number of those components $C$ of
%$G-S-T$ which is $H_{(X,Y)|C}$-critical.

\begin{theorem}\label{thm2}
Let $F_{i}=F[V(D_{i})]$ for $i=1,\ldots,\tau$. Then
$def_{H_{(A,B)|D_{i}}}(F_{i})=1$ and $F_{i}$ is
$H_{(A,B)|D_{i}}$-optimal for $i=1,\ldots,\tau$.
\end{theorem}
\vspace{-3mm} \pf By Lemmas \ref{prop1} and \ref{prop4}, we have
$def_{H_{(A,B)|D_{i}}}(F_{i})\leq 1$. Suppose that the theorem
doesn't hold. Let $F_{i}^*$ be $H_{(A,B)|D_{i}}$-optimal. Then we
have $def_{F_{i}^*}(D_{i})=0$.
 We consider two cases.

\medskip
{\it Case 1.} $def_{H}(F;D_{i})=1$.

Then we have $def_{H}((F-F_{i})\cup F_{i}^*)<def_{H}(F)$, but $F$ is $H$-optimal, a contradiction.

\medskip
{\it Case 2.} $F$ contains (or misses) an edge of $E(D_{i},A)$
(resp. $E(D_{i},B)$).

Let $yy_{1}\in E(D_{i},A)\cap E(F)$, where $y\in D_{i}$ and $y_{1}\in A$ (resp. $yy_{1}\in
E(D_{i},B)$ and $yy_{1}\notin E(F)$, where $y\in D_{i}$ and $y_{1}\in B$). Since $y\in D_{i}$,
there is a changeable trail $P$ from a vertex $x$ of $B_{0}$ to $y$.  By Lemma \ref{prop4}, we have
$yy_{1}\in E(P)$. Let $P_{1}$ be a subtrail of $P$ such that
 $V(P_{1})\cap V(D_{i})=\{y\}$. Then we have
$def_{H}(((F-F_{i})\triangle P_{1})\cup F_{i}^*)<def_{H}(F)$, a
contradiction again.
 \qed

Now we prove Lov\'asz's classic deficiency formula.

\begin{theorem}[Lov\'asz \cite{lov1}]\label{Lovasz}
The total deficiency is
$$def_H(G) =\max \limits_{S, T \atop S\cap T = \emptyset} \tau_H(S, T) - \sum_{x\in T}d_{G-S}(x) - MH(S) + mH(T),$$
where $\tau_H(S,T)$ denotes the number of components $K$ of $G-S-T$ such that $K$ contains no
$H_{(S,T)|K}$-factors. Moreover, a graph $G$ has an $H$-factor if and only if for any pair of
disjoint sets $S, T \subseteq V(G)$,
$$\tau_H(S, T) - \sum_{x\in T}d_{G-S}(x) - MH(S) + mH(T)\leq 0.$$
\end{theorem}

\pf Let $M$ be an arbitrary  $H$-optimal graph of $G$. Firstly, we
show that
$$def_H(G) \geq\max \limits_{S, T \atop S\cap T = \emptyset}
\tau_H(S, T) - \sum_{x\in T}d_{G-S}(x) - MH(S) + mH(T).$$ Let $S$
and $T$ be arbitrary disjoint subsets of $V(G)$. Let $\tau_{H}(S,T)$
be defined as in the above. For $i=1,,\ldots, \tau_{H}(S,T)$, let
$C_{i}$ denote the component of $G-S-T$ containing no
$H_{(S,T)|C_{i}}$-factors. Let $W=C_{1}\cup\cdots\cup
C_{\tau_{H}(S,T)}$. Since $C_{i}$ contains no
$H_{(S,T)|C_{i}}$-factors, so if $def_{M}(C_{i})=0$, then $M$ either
misses at least an edge of $E(C_{i},T)$ or contains at least an
edges of $E(C_{i},S)$. Let $\tau_{1}$ denote the number of
components of $W$ such that $M$ misses at least an edge of
$E(C_{i},T)$ and $\tau_{2}$ denote the number of the components of
$W$  such that $M$ contains at least an edge of $E(C_{i},S)$. Then
we have
\begin{align*}
def_{H}(G)&=def_{H}(M)\geq
\tau_{H}(S,T)-\tau_{1}-\tau_{2}+\sum_{x\in S\cup
T}dist(d_{M}(x),H(x))\\
&\geq \tau_{H}(S,T)-\tau_{1}-\tau_{2}+\sum_{x\in
S}(d_{M}(x)-MH(x))+\sum_{x\in
T}(mH(x)-d_{M}(x))\\
&\geq
\tau_{H}(S,T)-\tau_{1}-\tau_{2}+(e_{M}(S,T)+\tau_{2}-MH(S))+\sum_{x\in
T}(mH(x)-d_{M}(x))\\
&=\tau_{H}(S,T)-\tau_{1}+(e_{M}(S,T)-MH(S))+\sum_{x\in
T}(mH(x)-d_{M}(x))\\
&\geq
\tau_{H}(S,T)-\tau_{1}+(e_{M}(S,T)-MH(S))+(mH(T)-(e_{M}(S,T)+\sum_{x\in
T}d_{G-S}(x)-\tau_{1}))\\
&=\tau_{H}(S,T)+mH(T)-MH(S)-\sum_{x\in T}d_{G-S}(x).
\end{align*}
By Theorems \ref{thm1}, we have
$$def_{H}(G)=\tau+\sum_{v\in B}(mH(v)-d_{G-A}(v))-\sum_{v\in
A}MH(v).$$ By Theorem  \ref{thm2}, $D_{i}$  contains no
$H_{(A,B)|D_{i}}$-factors for $i=1,\cdots,\tau$. So
 we have
\begin{align*}
def_H(G) &=\tau_H(A, B)+\sum_{v\in B}(mH(v)-d_{G-A}(v))-\sum_{v\in
A}MH(v)\\
&= \max \limits_{S, T \atop S\cap T = \emptyset} \tau_H(S,
T) - \sum_{x\in T}d_{G-S}(x) - MH(S) + mH(T).
\end{align*}
 We complete the proof.
\qed

The proof of Theorem \ref{Lovasz} also imply the following result.

\begin{theorem}\label{LuYu}
Let $R$ be an arbitrary $H$-optimal graph of $G$. Then
\begin{itemize}
\item[$(1)$] $d_{R}(v)\in H(v)$ for all $v\in C$;

\item[$(2)$] $d_{R}(v)\geq MH(v)$ for all $v\in A$;

\item[$(3)$] $d_{R}(v)\leq  mH(v)$ for all $v\in B$.
\end{itemize}
\end{theorem}

From Theorem \ref{LuYu}, we can see that the partition $(A(G), B(G),
C(G), D(G))$ defined in this paper is equivalent to the original
partition $((A_H, B_H, C_H, D_H)$ introduced by Lov\'asz in
\cite{lov1}.

\begin{theorem}
$C_{H}=C$, $A_{H}=A$, $B_{H}=B$ and $D_{H}=D$.
\end{theorem}

\pf By Theorem \ref{LuYu}, we have $C\subseteq C_{H}$, $A\subseteq
A_{H}$, $B\subseteq B_{H}$. However, the definition of $D$ implies
$v\notin C_{H}\cup A_{H}\cup B_{H}$ for every $v\in D$. So we have
$D\subseteq D_{H}$. This completes the proof. \qed

\end{document}